\documentclass[12pt]{amsart}
\usepackage{amscd,amssymb,mathrsfs,epsf,subfigure,float}
\usepackage{xy}\xyoption{all}
\newtheorem*{thm*}{Theorem}
\newtheorem{thm}{Theorem}[section]
\newtheorem{conj}[thm]{Conjecture}

\newcommand{\C}{\mathbb C}
\newcommand{\Fix}{\mathsf{Fix}}
\newcommand{\R}{\mathbb R}
\newcommand{\Z}{\mathbb Z}
\newcommand{\EE}{\mathsf E}
\renewcommand{\L}{\mathsf L}
\newcommand{\V}{\mathsf V}
\newcommand{\tM}{\tilde M}
\newcommand{\dev}{\mathsf{dev}}
\newcommand{\Aff}{\mathsf{Aff}}
\newcommand{\GL}{\mathsf{GL}}
\newcommand{\GLp}{\GL^+(2,\R)}
\newcommand{\tGLp}{\widetilde{\mathsf{GL}^+(2,\R)}}
\newcommand{\rpo}{\R\mathsf{P}^1}
\newcommand{\SOt}{\mathsf{SO}(2)}
\newcommand{\Ht}{\mathsf{H}^2}
\newcommand{\Isom}{\mathsf{Isom}}
\newcommand{\PSLtR}{\mathsf{PSL}(2,\R)}
\newcommand{\Hom}{\mathsf{Hom}}

\newcommand{\tG}{\tilde{G}}

\newcommand{\Eto}{\mathbb{E}^3_1}
\newcommand{\Uoo}{\mathsf{U}(1,1)}
\newcommand{\Coo}{\mathfrak{u}(1,1)}
\newcommand{\PUto}{\mathsf{PU}(2,1)}
\newcommand{\PUoo}{\mathsf{PU}(1,1)}
\newcommand{\SUoo}{\mathsf{SU}(1,1)}
\newcommand{\Uo}{\mathsf{U}(1)}
\newcommand{\Ad}{\mathsf{Ad}}
\newcommand{\gar}{\mathfrak{g}_{\Ad\rho}}
\newcommand{\htc}{\mathsf{H}^3_\C}
\newcommand{\hoc}{\mathsf{H}^3_\C}
\newcommand{\SLtR}{\mathsf{SL}(2,\R)}

\newcommand{\Oo}{\mathsf{O}}
\newcommand{\oto}{\Oo(2,1)}
\newcommand{\CS}{\mathscr{C}(\Sigma)}
\newcommand{\ml}{\mathscr{ML}(\Sigma)}
\newcommand{\ho}{H^1(\Gamma_0,\V)}
\newcommand{\Proper}{\mathsf{Proper}}
\newcommand{\Fricke}{\mathfrak{F}(\Sigma)}
\newcommand{\W}{\mathsf{W}}
\newcommand{\Li}{L_\infty}
\newcommand{\gad}{\mathfrak{g}_{\Ad}}
 
\begin{document}

\title[Milnor's work on flat manifolds]
{Two papers which changed my life:\\  
Milnor's seminal work on flat manifolds
and bundles}
\author[W. Goldman]{William M. Goldman}
\date{\today}  

\thanks{
This paper was presented at the workshop ``Frontiers in 
Complex Dynamics'' at the Banff International Research Station,
in Banff, Alberta, Canada. 
I gratefully acknowledge partial support from National Science Foundation
grant DMS070781.}

\begin{abstract}
We survey developments arising from 
Milnor's 1958 paper, ``On the existence of a connection
with curvature zero'' and his 1977 paper, ``On fundamental groups of complete
affinely flat manifolds.''
\end{abstract}
\maketitle
\setcounter{tocdepth}{1} 
\dedicatory{\it With warm wishes to Jack Milnor on his eightieth birthday}
\tableofcontents

For a young student studying topology at Princeton in the
mid-1970's, John Milnor was a inspiring presence.  The excitement of
hearing him lecture at the Institute for Advanced Study and reading his
books and unpublished lecture notes available in Fine Library made a deep impact on me.  One heard rumors of exciting breakthroughs in the
Milnor-Thurston collaborations on invariants of $3$-manifolds
and the theory of kneading in $1$-dimensional dynamics.
The topological significance of volume in hyperbolic $3$-space
and Gromov's  proof of Mostow rigidity using simplicial
volume were in the air at the time (later to be written up in Thurston's
notes~\cite{Thurston}).
When I began studying 
geometric structures on manifolds,
my mentors Bill Thurston and Dennis Sullivan 
directed me to two of Milnor's papers~\cite{Milnor,Milnor_connection}.
Like many mathematicians of my generation, 
his papers and books such as {\em Morse theory,\/}  
{\em Characteristic Classes,\/} and 
{\em Singular Points of Complex Hypersurfaces,\/} 
were very influential for my training.
Furthermore his lucid writing style made his papers and books 
role models for exposition.

I first met Jack in person when I was a graduate student in Berkeley and he was visiting his son.  Several years later I was extremely flattered when I received a letter from him,
where he, very politely, pointed out a technical error in my Bulletin Announcement~\cite{GoldmanBAMS}.
It is a great pleasure and honor for me to express my  gratitude
to Jack Milnor for his inspiration and insight, in celebration of his
eightieth birthday.

\section{Gauss-Bonnet beginnings}

\subsection{Connections and characteristic classes}
The basic topological invariant of a closed orientable surface $M^2$
is its Euler characteristic $\chi(M)$. If $M$ has genus $g$, then
$\chi(M) = 2- 2g$. 
Give $M$ a Riemannian metric --- then
the Gauss-Bonnet theorem identifies $\chi(M^2)$ as $1/2\pi$ of
the total Gaussian curvature of $M$, a geometric invariant.
This provides a fundamental topological restriction of what kind of geometry $M$ may support.

For example, if the metric is flat --- that is, locally Euclidean --- then
its Gaussian curvature vanishes and therefore $\chi(M) = 0$. 
Since $M$ is orientable,  $M$ must be homeomorphic to a torus.

In 1944 Chern~\cite{Chern} proved his intrinsic Gauss-Bonnet theorem,
which expresses the topological invariant $\chi(M^n)$ as the integral
of a differential form constructed from the curvature of a Riemannian
metric. More generally, if $\xi$ is an oriented $n$-plane bundle 
over $M^n$, then its {\em Euler number\/} 
$$
e(\xi)\in H^n(M,\Z)\cong \Z
$$
can be computed as an integral of an expression derived from
an {\em orthogonal connection\/} $\nabla$ on $\xi$.
(Compare Milnor-Stasheff~\cite{CharacteristicClasses}.) 
For example take to $M$ to be a (pseudo-) Riemannian manifold, with 
tangent bundle $\xi = TM$ and $\nabla$ the Levi-Civita connection.
If $\nabla$ has curvature zero, then, according to Chern,
$\chi(M) = 0$.

This paradigm generalizes.
When $M$ is a complex manifold, its tangent bundle
is a holomorphic vector bundle and the Chern classes can
be computed from the curvature of 
a holomorphic connection. In particular $\chi(M)$
is a Chern number. Therefore, if $TM$ has a flat holomorphic
connection, then $\chi(M) = 0$. If $M$ has a flat pseudo-Riemannian
metric, then a similar Gauss-Bonnet theorem holds (Chern~\cite{Chern2}) and $\chi(M)=0$. 
Therefore it is natural to ask whether a compact manifold whose tangent bundle admits a flat  
{\em linear connection\/} has Euler characteristic zero.

\subsection{Smillie's examples}
In 1976, John Smillie~\cite{Smillie} constructed, in every even
dimension $n > 2$, a compact $n$-manifold such that $TM$ admits
a flat connection $\nabla$. However, the torsion of $\nabla$
is (presumably) nonzero. 
(None of Smillie's examples are {\em aspherical;\/}
it would be interesting to construct a closed aspherical manifold
with flat tangent bundle; compare Bucher-Gelander~\cite{BucherGelander}.)

Requiring the torsion of $\nabla$ to vanish is a natural condition.
When both curvature and torsion vanish, the connection arises
from an {\em affine structure\/} on $M$, that is, the structure defined
by a coordinate atlas of coordinate charts into an
affine space $\EE$ such that the coordinate changes on overlapping
coordinate patches are locally affine. A manifold together
with such a geometric structure is called an {\em affine manifold.\/}  
The coordinate charts globalize into a {\em developing map\/}
$$
\tM \xrightarrow{\dev} \EE
$$
where $\tM\longrightarrow M$ is a universal covering space.
The developing map
$\dev$ is a local diffeomorphism (although generally not a covering
space onto its image), which defines the affine structure.
Furthermore $\dev$ is equivariant with respect to a homomorphism
$$
\pi_1(M) \xrightarrow{\rho} \Aff(\R^2)
$$
(the {\em affine holonomy representation\/}) where
the fundamental group $\pi_1(M)$ acts by deck transformations
of $\tM$. Just as $\dev$ globalizes the coordinate charts,
$\rho$ globalizes the {\em coordinate changes.\/}
The flat connection on $TM$ arises from the representation
$\rho$ in the standard way: 
$TM$ identifies with the fiber product 
$$
\xi_\rho \;:=\;  \tM \times_\rho \R^n \; = \; (\tM \times \R^n)/\big(\pi_1(M)\big)
$$
where $\pi_1(M)$ acts {\em diagonally\/} ---  by deck transformations
on the $\tM$ factor and via $\rho$ on the $\R^n$-factor.
The differential of $\dev$ defines an isomorphism of $TM$ with 
$\xi_\rho$.

We may interpret Smillie's examples in this description as follows.
In each even dimension $n>2$, Smillie constructs an $n$-manifold
$M^n$ and a representation $\pi_1(M)\xrightarrow{\rho}\Aff(\R^n)$
such that $\xi_\rho$ is isomorphic to the tangent bundle $TM$.
Sections $\xi_\rho$ correspond to {\em singular developing maps\/}
which may be smooth, but {\em not\/}%
necessarily %
local diffeomorphisms.

Despite %
the %
many partial results, we know
no example of a closed affine manifold with 
{\em nonzero Euler characteristic.}
%

\subsection{Benz\'ecri's theorem on flat surfaces}
In dimension two, a complete answer is known, 
due to the work of Benz\'ecri~\cite{Benzecri}.
This work was part of his 1955 thesis at Princeton,
and Milnor served on his thesis committee.

\begin{thm*}[Benz\'ecri]
A closed surface $M$ admits
an affine structure if and only if $\chi(M)=0$. 
\end{thm*}
\noindent
(Since every connected orientable open surface can be immersed
in $\R^2$, pulling back the affine structure from $\R^2$ by
this immersion gives an affine structure. With a small modification
of this technique, every connected nonorientable open surface can
also be given an affine structure.)

Benz\'ecri's proof is geometric, and starts with a fundamental
polygon $\Delta$ for $\pi_1(M)$ acting on $\tM$. 
The boundary $\partial\Delta$ consists of various edges, which
are paired by homeomorphisms, reconstructing $M$ as the quotient
space by these identifications.
One standard setup for a surface of genus $g$ uses a $4g$-gon
for $\Delta$, where the sides are alternately paired to give
the presentation 
\begin{equation}\label{eq:presentation}
\pi_1(M) \;=\; \langle A_1, B_1, \dots, A_g, B_g \mid
A_1 B_1A_1^{-1} B_1^{-1} \dots A_g B_g A_g^{-1}B_g^{-1} = 1 \rangle
\end{equation}
The developing map $\dev$ immerses $\Delta$ into $\R^2$, and
the identifications between the edges of $\partial\Delta$ are realized
by orientation-preserving affine transformations. 

Immersions of $S^1$ into $\R^2$ are classified up to regular homotopy
by their {\em turning number\/} 
(the Whitney-Graustein theorem\cite{Whitney}) 
which measures the total angle the tangent vector (the velocity) 
turns as the curve is traversed. 
Since the restriction $\dev|_{\partial\Delta}$ extends to an immersion of the disc $\Delta$ its turning
number (after choosing compatible orientations) 
\begin{equation}\label{tauBounds}
\tau(\dev|_{\partial\Delta}) = 2\pi.
\end{equation}
However, Benz\'ecri shows that for any smooth immersion
$[0,1]\xrightarrow{f}\R^2$ and orientation-preserving affine
transformation $\gamma$,
\begin{equation}\label{eq:tauLess}
\vert \tau(f) - \tau(\gamma\circ f) \vert < \pi.
\end{equation}
Using the fact that $\partial\Delta$ consists of $2g$ pairs of edges
which are paired by $2g$ orientation-preserving affine transformations,
combining \eqref{tauBounds} and \eqref{eq:tauLess} implies $g=1$.

Milnor realized the algebraic-topological ideas underlying 
Benz\'ecri's proof, thereby initiating the theory of characteristic classes
of flat bundles.

\section{The Milnor-Wood inequality}

\subsection{``On the existence of a connection with curvature zero''}
In his 1958 paper~\cite{Milnor}, Milnor shows that a closed
$2$-manifold $M$ has flat tangent bundle if and only if
$\chi(M) =0$. This immediately implies Benz\'ecri's theorem, 
although it doesn't use the fact that the developing map
is nonsingular (or, equivalently, the associated flat connection 
is torsionfree). 
In this investigation, Milnor discovered, remarkably,
flat oriented $\R^2$-bundles over surfaces $M$ with 
{\em nonzero Euler class.\/}  
In particular the Euler class 
{\em cannot\/} be computed from the curvature of a linear connection.

Oriented $\R^2$-bundles $\xi$ over $M$ are
classified up to isomorphism by their Euler class 
$$
e(\xi) \in H^2(M;\Z)
$$
and if $M$ is an orientable surface, an orientation on $M$
identifies $H^2(M;\Z)$ with $\Z$. 
(See Milnor-Stasheff~\cite{CharacteristicClasses} for details.)
An oriented $\R^2$-bundle $\xi$ admits a flat structure if and only
if it arises from a  representation 
$$
\pi_1(M) \xrightarrow{\rho} \GLp
$$
(where $\GLp$ denotes the group of orientation-preserving
linear automorphisms of $\R^2$).
Milnor shows that $\xi$ admits a flat structure 
if and only if its Euler number satisfies
\begin{equation}\label{eq:MilnorIneq}
\vert e(\xi) \vert < g.
\end{equation}
Since $e(TM) = \chi(M) =  2-2g$, 
Milnor's inequality \eqref{eq:MilnorIneq} implies that $g = 1$.

The classification of $S^1$-bundles is basically equivalent to the
classification of rank $2$ vector bundles, but is somewhat more
general.  To any vector bundle $\xi$ with fiber $\R^n$ is associated
an $S^{n-1}$-bundle: the fiber of the $S^{n-1}$-bundle over a point
$x$ consists of all directions in the fiber $\xi_x\approx \R^n$.  In
particular two $\R^2$-bundles are isomorphic if and only if their
associated $S^1$-bundles are isomorphic.  Therefore we henceforth work
with $S^1$-bundles, slightly abusing notation by writing $\xi$ for the
$S^1$-bundle associated to $\xi$.

\subsection{Wood's extension and foliations}
In 1971, John W.\ Wood~\cite{Wood} extended Milnor's classification
of flat $2$-plane bundles to flat $S^1$-bundles.
Circle bundles with structure group $\GLp$ have an important
special property. 
The {\em antipodal map\/} associates a direction in a vector
space its opposite direction.
Since all linear transformations commute with it,
the antipodal map defines an involution on any vector bundle
or associated sphere bundle $\xi$.
The quotient is the associated $\rpo$-bundle $\hat\xi$,
and the quotient map $\xi \longrightarrow \hat\xi$ is a double covering. 
This is also an oriented $S^1$-bundle, with Euler class
$$
e(\hat\xi) = 2 e(\xi).
$$
Wood~\cite{Wood} determines the flat oriented $S^1$-bundles,
for an arbitrary homomorphism
$$
\pi_1(\Sigma) \xrightarrow{\rho} \mathsf{Homeo}^+(S^1).
$$
He proves the Euler number satisfies
the following inequality:
\begin{equation}\label{eq:MilnorWood}
\vert e(\rho) \vert \le -\chi(\Sigma),
\end{equation}
(now known as the {\em Milnor-Wood inequality\/}). %
Furthermore every integer in $[\chi(\Sigma),-\chi(\Sigma)]$
occurs as $e(\rho)$ for some homomorphism $\rho$.

Milnor's proof interprets the Euler class as the obstruction for lifting
the holonomy representation $\rho$ from the group $\GLp$
of linear transformations of $\R^2$ with positive determinant
to its universal covering group $\tGLp$.
Suppose $G$ is a Lie group with universal covering
$\tG \longrightarrow G$.
If $S_g$ is a closed oriented surface of genus $g>1$,
then its fundamental group admits a presentation
\eqref{eq:presentation}
Let 
$$
\pi_1(S_g) \xrightarrow{\rho} G
$$ 
be a representation; then the {\em obstruction\/}
$o_2(\rho)$ for lifting $\rho$ to $\tG$ is 
obtained as follows.
Choose lifts $\widetilde{\rho(A_i)}$ ,
$\widetilde{\rho(B_i)}$ of $\rho(A_i)$ and $\rho(B_i)$ 
to $\tG$ respectively.
Then
\begin{equation}\label{eq:CommProd}
o_2(\rho) := [\widetilde{\rho(A_1)},\widetilde{\rho(B_1)}]\dots
[\widetilde{\rho(A_g)},\widetilde{\rho(B_g)}]
\end{equation}
is independent of the chosen lifts, and lies in 
$$
\pi_1(G) = \ker(\tG \longrightarrow G).
$$
It vanishes precisely when $\rho$ lifts to $\tG$.
When $G = \GLp$, the obstruction class $o_2(\rho)$
is just the Euler class $e(\rho)$. 

To identify the element of
$\pi_1(\GLp)$ corresponding to $e(\rho)$,
Milnor and Wood estimate the {\em translation number\/}
of the lifts of generators to $\tG$, which is based on the
rotation number of orientation-preserving circle homeomorphisms.
Milnor uses a retraction $\GLp\xrightarrow{r}\SOt$
(say, the one arising from the Iwasawa decomposition), 
which lifts to a retraction
$$
\tG \xrightarrow\theta \widetilde{\SOt} \cong \R
$$
and proves the estimate
\begin{equation}\label{eq:MilnorEstimate}
\vert \theta(\gamma_1\gamma_2) - \theta(\gamma_1) - \theta(\gamma_2) \vert < \frac{\pi}2.
\end{equation}
Wood considers a more general retraction 
$\theta$ defined on $\tG = \widetilde{\mathsf{Homeo}^+(S^1)}$
and shows a similar estimate, sharpened by a factor of two.
Applying this to 
\eqref{eq:CommProd}, he shows that
if an $m$-fold product of commutators in $\tG$
is translation by $a$, then
\begin{equation}\label{eq:WoodEstimate}
\vert a\vert < 2 m - 1
\end{equation}
The estimate \eqref{eq:WoodEstimate}
extends  Benz\'ecri's original estimate \eqref{eq:tauLess}
in a stronger and  more abstract context. 
This --- the boundedness of the Euler class of flat bundles ---
may be regarded as one of the roots of the theory
of bounded cohomology. 
The fundamental role of the Euler class as a {\em bounded
cohomology class\/} was discovered by Ghys~\cite{Ghys}.
In particular he showed that the {\em bounded Euler class\/}
characterizes orientation-preserving
actions of surface groups on the circle up to quasi-conjugacy.

For other generalizations of 
the Milnor-Wood inequality, compare
Dupont~\cite{Dupont1,Dupont2},
Sullivan~\cite{Sullivan}, Domic-Toledo~\cite{DomicToledo} 
and Smillie~\cite{Smillie_unpublished}.
For more information, see 
Burger-Iozzi-Wienhard~\cite{BurgerIozziWienhard_handbook}
and  the second chapter of Calegari~\cite{Calegari}.
The question of when a foliation on the total space of a circle bundle over a surface is isotopic to  a flat bundle is the subject of 
Thurston's thesis~\cite{Thurston_thesis}.

\section{Maximal representations}

\subsection{A converse to the Milnor-Wood inequality}

Equality in \eqref{eq:MilnorWood} has special and deep
significance. 
Let $M$ be a closed oriented surface.
Then, just as  described earlier for affine structures, 
every hyperbolic structure on $M$ determines a 
{\em developing pair\/} $(\dev,\rho)$ where
\begin{align*}
\tM &\xrightarrow{\dev} \Ht \\
\pi_1(M) &\xrightarrow{\rho} \Isom^+(\Ht).
\end{align*}
by globalizing the coordinate charts and coordinate charts
in an atlas defining the hyperbolic structure.
The flat $(\Isom^+(\Ht),\Ht)$-bundle $E_M \longrightarrow M$
corresponding to $\rho$ has a section $\delta_M$ corresponding to 
$\dev$,
which is transverse to the flat structure on $E_M$. 
Consequently the normal bundle of $\delta_M\subset E_M$
(by the tubular neighborhood theorem) is isomorphic to the
tangent bundle $TM$, and therefore 
$$
e(\rho) = e(TM) = \chi(M),
$$
proving sharpness in the Milnor-Wood inequality.
By conjugating $\rho$ with an orientation-reversing isometry of $\Ht$,
one obtains a representation $\rho$ with 
$e(\rho) = - \chi(M)$.

The converse statement was proved in my doctoral 
dissertation~\cite{Goldman_thesis}.
Say that a representation is {\em maximal\/}
if $e(\rho) = \pm \chi(M)$. 

\begin{thm*}
Let $\rho\in\Hom(\pi_1(M),\PSLtR)$. 
Then the following are equivalent:
\begin{itemize}
\item $\rho$ is the holonomy of a hyperbolic structure on $M$;
\item $\rho$ is an embedding onto a discrete subgroup of 
$\Isom^+(\Ht)\cong \PSLtR$;
\item For every $\gamma\in\pi_1(M)$ with $\gamma\neq 1$,
the holonomy $\rho(\gamma)$ is a hyperbolic element of 
$\PSLtR$.
\end{itemize}
\end{thm*}

\subsection{Kneser's theorem on surface maps}

A special case follows from the classical theorem of Kneser~\cite{Kneser}:

\begin{thm*}
Let $M,N$ be closed oriented surfaces, and $N$ having genus $>1$.
Suppose that $M\xrightarrow{f} N$ is a continuous map of degree $d$.
Then
$$
d \vert\chi(N)\vert \le \vert\chi(M)\vert.
$$
Furthermore
$d \vert\chi(N)\vert = \vert\chi(M)\vert$
if and only if $f$ is homotopic to a covering space. %
\end{thm*}
\noindent 
The theorem follows by giving $N$ a hyperbolic structure,
with holonomy representation $\rho$. 
Then the composition
$$
\pi_1(M) \xrightarrow{f_*} \pi_1(N) \xrightarrow{\rho}  \PSLtR
$$
has Euler number $d\chi(N)$. Now apply Milnor-Wood and its
converse statement to the composition.

\subsection{Components of the representation variety}
Since the space of hyperbolic structures on $M$ is connected,
the Euler class defines a continuous map
$$
\Hom\big(\pi_1(M),\PSLtR\big) \longrightarrow H^2(M;\Z) \cong \Z.
$$
Reversing the orientation on $M$ reverses the sign of 
the Euler number. 
Therefore the maximal representations constitute {\em two}
connected components of $\Hom\big(\pi_1(M),\PSLtR\big)$.
In general the connected components are the fibers of this
map (Goldman~\cite{GoldmanBAMS,GoldmanINV},
Hitchin~\cite{Hitchin}).
In particular the space of representations has $4g-3$ connected
components.
Each component has dimension $6g-6$.
Only the Euler class $0$ component is not a smooth manifold.
Furthermore Hitchin relates the component corresponding to
Euler class $2-2g + k$ to the $k$-th symmetric power of $M$.
See Hitchin~\cite{Hitchin},
Bradlow-Garc\`{\i}a-Prada- Gothen~\cite{BradlowGarciaPradaGothen2,
Bradlow_GarciaPrada_Gothen_whatis}, 
as well as my  expository article \cite{Goldman_Hitchin}.

\subsection{Rigidity and flexibility}

This characterization of maximal representations is a kind of rigidity
for surface group representations.  Dupont~\cite{Dupont1,Dupont2},
Turaev~\cite{Turaev,Turaev1} and Toledo~\cite{Toledo} defined
obstruction classes $o_2$ for Lie groups $G$ of automorphisms of
Hermitian symmetric spaces.  In particular Toledo~\cite{Toledo} proved
the following rigidity theorem:
\begin{thm*}{\rm(Toledo~\cite{Toledo})}
Suppose that $\pi_1(M)\xrightarrow{\rho}\mathsf{PU}(n,1)$ is a representation.
Equality is attained in the generalized Milnor-Wood
inequality
$$
\vert o_2(\rho) \vert \le \frac{\vert\chi(M)\vert}2.
$$
Then $\rho$ embeds $\pi_1(M)$ as a discrete subgroup of the stabilizer
(conjugate to $\mathsf{U}(1,1)\times \mathsf{U}(n-1)$ of a holomorphic 
totally geodesic curve $C$ in $\mathbf{H}_\C^2$.
In particular $C/\mathsf{Image}(\rho)$ is a hyperbolic surface
diffeomorphic to $M$.
\end{thm*}
\noindent
Recently these results have been extended to higher rank
in the work of  Burger-Iozzi-Wienhard~\cite{Burger_Iozzi_Labourie_Wienhard,
BurgerIozziWienhard_toledo,BurgerIozziWienhard_handbook}.

As maximality of the Euler class in the Milnor-Wood inequality
implies rigidity, various values of the Euler class imply various
kinds of flexibility~\cite{GoldmanRep}. 
If $\pi = \pi_1(M)$ is the fundamental group of a compact 
K\"ahler manifold $M$, and $G$ is a reductive algebraic Lie group,
then Goldman-Millson~\cite{GoldmanMillson} gives a complete
description of the the analytic germ of the space of representations 
$\Hom(\pi,G)$
at a reductive representation $\rho$.
Specifically, $\rho$ has an open neighborhood in $\Hom(\pi,G)$
analytically equivalent to the quadratic cone
defined by the symmetric bilinear form
$$
Z^1(\pi,\mathfrak{g}_{\mathsf{Ad}\rho}) \times
Z^1(\pi,\mathfrak{g}_{\mathsf{Ad}\rho}) \longrightarrow
H^2(\pi,\mathfrak{g}_{\mathsf{Ad}\rho}) 
$$
obtained by combining cup product on $\pi$ with Lie bracket
$$
\mathfrak{g}_{\mathsf{Ad}\rho}\times
\mathfrak{g}_{\mathsf{Ad}\rho} \xrightarrow{[,]}
\gar
$$
as coefficient pairing.

Consider the special case when $M$ is a closed hyperbolic surface
and  a representation $\pi \xrightarrow{\rho_0} 
\SUoo\cong\SLtR$.
We assume that $\rho_0$ has Zariski-dense image, 
which in this case simply means that its image is non-solvable.
In turn this means %
the corresponding action on $\Ht$ fixes no point in $\Ht\cup\partial\Ht$.
In that case $\rho_0$ is reductive and defines a smooth point
of the $\R$-algebraic set  $\Hom(\pi,\SUoo)$. %
Extend the action to an isometric action on the complex hyperbolic
plane $\htc$ via the composition $\rho$ defined by:
\begin{equation}\label{eq:composite}
\pi \xrightarrow{\rho_0}  \mathsf{SU}(1,1) \hookrightarrow
\PUto
\end{equation}
For $A\in \mathsf{U}(1,1)$, 
taking the equivalence class of the direct sum
$$
A \oplus 1 := \bmatrix A & 0 \\ 0 & 1 \endbmatrix
\in \mathsf{U}(2,1) $$
in $\PUto$  defines an  embedding 
$$
\mathsf{SU}(1,1) \hookrightarrow\PUto.
$$
This representation stabilizes a {\em complex hyperbolic line\/}
$\hoc \subset \htc$ inside the complex hyperbolic plane.
What are the local deformations of $\rho$ in
$\Hom(\pi,\PUto)$?

The representation $\rho$ is maximal if and only if $\rho_0$
is maximal, which occurs when $e(\rho) = \pm \chi(M)$.
In that case any representation $\pi\longrightarrow\PUto$
near $\rho$ stabilizes $\hoc$, that is, it lies in the subgroup
$\Uoo\subset\PUto$. 

In general, representations $\pi\longrightarrow\Uoo$
can be easily understood in terms of their composition with the
projectivization %
homomorphism $\Uoo\longrightarrow\PUoo$.
The corresponding map on representation varieties
$$
\Hom(\pi,\Uoo) \longrightarrow \Hom(\pi,\PUoo)
$$ 
is a torus fibration, where the points of the fiber correspond to
different actions in the normal directions to $\hoc\subset\htc$
(which are described by characters $\pi\longrightarrow \Uo$.
In particular $\rho$ defines a smooth point of $\Hom(\pi,\Uoo)$ 
with tangent space 
$$
Z^1\big(\pi,\mathfrak{su}(1,1)_{\Ad\rho_0} \big) \oplus
Z^1\big(\pi,\R)
$$
since $\mathfrak{u}(1)_{\Ad\rho_0})$ equals the ordinary
coefficient system $\R$ (where $\pi$ acts by the identity).

The general deformation result implies that the analytic germ
of $\Hom(\pi,\PUto)$ near $\rho$ looks like the Cartesian product
of the (smooth) analytic germ of 
$$
\Hom(\pi,\Uoo) \times \PUto/\Uoo
$$ 
with the quadratic cone 
$\mathscr{Q}_{\rho}$ in $Z^1(\pi,\Coo_{\Ad \rho}) \cong \C^{2g}$
defined by the cup-product pairing
$$
Z^1(\pi,\Coo_{\Ad \rho}) \times
Z^1(\pi,\Coo_{\Ad \rho}) \longrightarrow
H^2(\pi,\R) \cong\R
$$
where 
$\Coo_{\rho}$ denotes the $\pi$-module defined by the standard
$2$-dimen\-sional complex representation of $\Uoo$.
The coefficient pairing (which is derived from the Lie bracket
on $\mathfrak{su}(2,1)$) is just the imaginary part of the indefinite
Hermitian form on $\Coo$, and is skew-symmetric.
In particular the space of {\em coboundaries\/} 
$$
B^1(\pi,\Coo_{\Ad \rho})\subset Z^1(\pi,\Coo_{\Ad \rho}),
$$
is isotropic. Thus we reduce to  the symmetric
bilinear form obtained from the cohomology pairing

\begin{equation}\label{eq:CohomologyPairing}
H^1(\pi,\Coo_{\Ad\rho}) \times
H^1(\pi,\Coo_{\Ad\rho}) \longrightarrow H^2(\pi,\R) \cong\R.
\end{equation}
\noindent
The real dimension of 
$H^1(\pi,\Coo_{\rho})$ equals 
$$
-2\chi(M) = 8(g-1).
$$
By the Signature Theorem  of Meyer~\cite{Meyer}, 
the quadratic form corresponding to 
\eqref{eq:CohomologyPairing}
has signature $8 e(\rho_0)$.
Thus, near $\rho$, the $\R$-algebraic
set $\Hom(\pi,\PUto)$ is analytically equivalent to 
the Cartesian product of a manifold 
with a cone on $\R^{8(g-1)}$ defined by a quadratic form
of signature $8e(\rho)$.

Meyer's theorem immediately gives a proof of Milnor's inequality
\eqref{eq:MilnorIneq}, since the signature of a quadratic form is bounded
by the dimension of the ambient vector space.
Furthermore $\rho$ is maximal if and only if the quadratic form
is definite, in which the quadratic cone has no real points,
and any small deformation of $\rho$ must stabilize a complex geodesic.

\section{Complete affine manifolds}
We return to the subject of flat affine manifolds,
and the second~\cite{Milnor} of Milnor's papers on this subject.

\subsection{The Auslander-Milnor question}
An affine manifold $M$ is {\em complete\/}
if some (and hence every) developing map
is bijective. In that case $\tM$ identifies
with $\EE$, %
 and $M$ arises as the quotient
$\Gamma\backslash\EE$ by a discrete subgroup 
$\Gamma\subset\Aff(\EE)$ acting properly
and freely on $\EE$. The affine holonomy
representation 
$$
\pi_1(M) \stackrel{\rho}\hookrightarrow \Aff(\EE)
$$
embeds $\pi_1(M)$ onto $\Gamma$.

Equivalently, an affine manifold is complete
if and only if the corresponding affine
connection is {\em geodesically complete,\/}
that is, every geodesic extends infinitely 
in both directions.

A simple example of an {\em incomplete\/} affine structure
on a closed manifold is a {\em Hopf manifold $M$,\/}
obtained as the quotient of $\R^n\setminus\{0\}$ by a cyclic group 
$\langle A\rangle$.
Here the generator $A$ must be a {\em a linear expansion,\/} 
that is, an element $A\in \GL(n,\R)$ such that every eigenvalue has modulus $> 1$.
Such a quotient is diffeomorphic to $S^{n-1} \times S^1$.
A geodesic aimed at the origin winds seemingly faster and faster around
the $S^1$-factor, although it's travelling with zero acceleration
with respect to the flat affine connection. In finite time, it
``runs off the edge" of the manifold.

If $M = \Gamma\backslash\EE$ is a complete affine manifold, 
then $\Gamma\subset\Aff(\EE)$ is a discrete subgroup 
acting properly and freely on $\EE$. 
However, in the example above, $\langle A\rangle$ 
is a discrete subgroup which doesn't act properly.
A proper action of a discrete group is the usual notion of
a {\em properly discontinuous action.\/}  If the action is also
free (that is, no fixed points), then the quotient is a (Hausdorff)
smooth manifold, and the quotient map $\EE \longrightarrow \Gamma\backslash\EE$ 
is a covering space.  %
A properly discontinuous action whose quotient is compact
as well as Hausdorff is said to be {\em crystallographic,\/}
in analogy with the classical notion of a {\em crystallographic group:\/}
A {\em Euclidean crystallographic group\/} is a discrete cocompact
group of Euclidean isometries. Its quotient space is a Euclidean
orbifold. Since such groups act isometrically on metric spaces,
discreteness here does imply properness; this dramatically fails
for more general discrete groups of {\em affine transformations.}

L.\ Auslander~\cite{Auslander1} claimed to prove that
the Euler characteristic vanishes for a compact complete
affine manifold, but his proof was flawed. 
It rested upon the following question, 
which in \cite{FriedGoldman}, was demoted to a
``conjecture,'' and is now known as the 
``Auslander Conjecture'':

\begin{conj}\label{conj:AuslanderConjecture}
Let $M$ be a compact complete affine manifold.
Then $\pi_1(M)$ is virtually polycyclic.
\end{conj}
\noindent
In that case the affine holonomy group $\Gamma\cong\pi_1(M)$
embeds in a closed Lie subgroup $G\subset \Aff(\EE)$ satisfying:
\begin{itemize}
\item $G$ has finitely many connected components;
\item The identity component $G^0$ acts simply transitively on $\EE$.
\end{itemize}
Then $M = \Gamma\backslash \EE$ admits a finite covering space  %
$M^0 := \Gamma^0\backslash\EE$
where 
$$
\Gamma^0 := \Gamma \cap G^0.
$$ 
The simply transitive action of $G^0$ define a complete {\em left-invariant affine structure\/}
on $G^0$. (The developing map is just the evaluation map of this action.)  
Necessarily $G^0$ is a $1$-connected solvable Lie group and
$M^0$ is affinely isomorphic to the {\em complete affine solvmanifold\/} $\Gamma^0\backslash G^0$.
In particular $\chi(M^0) = 0$ and thus $\chi(M) = 0$.

This theorem is the natural extension of Bieberbach's theorems
describing the structure of flat Riemannian (or Euclidean) manifolds;
see Milnor~\cite{Milnor_Hilbert} for an exposition of this theory
and its historical importance. 
Every flat Riemannian manifold is finitely covered by a {\em flat torus,\/}
the quotient of $\EE$ by a lattice of translations. In the more general case,  $G^0$ plays the role of the group of translations of an affine space and the solvmanifold $M^0$ plays the role of the flat torus.
The importance of Conjecture~\ref{conj:AuslanderConjecture} is that it would provide a detailed and computable structure theory for compact complete affine manifolds.

Conjecture~\ref{conj:AuslanderConjecture} was established
in dimension $3$ in Fried-Goldman~\cite{FriedGoldman}.
The proof involves classifying the possible Zariski closures
$A\big(\L(\Gamma)\big)$ of the linear holonomy group inside
$\GL(\EE)$.
Goldman-Kamishima~\cite{GoldmanKamishima} prove 
Conjecture~\ref{conj:AuslanderConjecture}
for flat Lorentz manifolds.
Grunewald-Mar\-gu\-lis~\cite{GrunewaldMargulis} establish
Conjecture~\ref{conj:AuslanderConjecture}
when the Levi component
of $\L(\Gamma)$ lies in a real rank-one subgroup of $\GL(\EE)$.
See Tomanov~\cite{Tomanov1,Tomanov2} and
Abels-Margulis-Soifer~\cite{AbelsMargulisSoifer,AbelsMargulisSoifer2,
AbelsMargulisSoifer3} for further results. The conjecture is now  known in all dimensions
$\le 6$ (Abels-Margulis-Soifer~\cite{AbelsMargulisSoifer4}).

\subsection{The Kostant-Sullivan theorem}

Although Conjecture~\ref{conj:AuslanderConjecture} 
remains unknown in general, 
the question which motivated it was proved by
Kostant and Sullivan~\cite{KostantSullivan}.

\begin{thm*}[Kostant-Sullivan]
Let $M$ be a compact complete affine manifold.
Then $\chi(M) = 0$ .
\end{thm*}

\noindent
Their ingenious proof uses an elementary fact about free affine actions
and Chern-Weil theory. 
The first step is that 
if $\Gamma\subset\Aff(\EE)$ is a group of affine transformations
acting freely on $\EE$, then the Zariski closure $A(\Gamma)$
of $\Gamma$ in $\Aff(\EE)$ has the property that every element
$g\in A(\Gamma)$ has $1$ as an eigenvalue.
To this end 
suppose that $\Gamma\subset\Aff(\EE)$ acts freely.
Then solving for a fixed point 
$$ 
\gamma(x) = \L(\gamma) (x) + u(\gamma) = x
$$ 
implies that $\L(\gamma)$ has $1$
as an eigenvalue for every $\gamma\in\Gamma$.
Thus  every element $\gamma\in \Gamma$ satisfies the polynomial
condition
$$
\det( \L(\gamma) - I ) = 0
$$
which extends to the Zariski closure $A(\Gamma)$ 
of $\Gamma$ in $\Aff(\EE)$.

Next one finds a Riemannian metric (or more accurately, an orthogonal connection)
to which Chern-Weil applies.
Passing to a finite covering, using the finiteness of $\pi_0\big(A(\Gamma)\big)$,
we may assume the holonomy group lies in the identity component $A(\Gamma)^0$, 
which is a connected Lie group.
Since every connected Lie group deformation retracts to a maximal compact subgroup,
the structure group of $TM$ reduces from $\Gamma$ to a maximal compact subgroup $K\subset A(\Gamma)^0$.
This reduction of structure group gives an orthogonal connection 
$\nabla$ taking values in the Lie algebra $\mathfrak{k}$ of $K$. 
Since every compact group of affine transformations fixes a point, 
we may assume that $K\subset \GL(\EE)$.
Since every element of $A(\Gamma)$
(and hence $K$) has $1$ as an eigenvalue, 
every element of $\mathfrak{k}$ has determinant zero. 
Thus the Pfaffian polynomial (the square root of the determinant)
vanishes on  $\mathfrak{k}$. Since the curvature of $\nabla$
takes values in $\mathfrak{k}$,
and the Euler form is the Pfaffian of the curvature tensor, the
Euler form is zero. Now apply the Chern-Gauss-Bonnet theorem~\cite{Chern}.
Integrating over $M$ gives $\chi(M) = 0$, as claimed.

\subsection{``On fundamental groups of complete affinely flat manifolds'' }

In his 1977 paper~\cite{Milnor},  
Milnor set the record straight caused by the confusion surrounding Auslander's flawed proof of  Conjecture~\ref{conj:AuslanderConjecture}.
Influenced by Tits's work~\cite{Tits} on free subgroups of linear groups
and amenability, Milnor observed, that for an affine space $\EE$ of given dimension, the following conditions are all equivalent:

\begin{itemize}
\item
Every discrete subgroup of $\Aff(\EE)$ which acts properly on $\EE$
is  amenable.
\item
Every discrete subgroup of $\Aff(\EE)$ which acts properly on $\EE$
is virtually solvable. 
\item
Every discrete subgroup of $\Aff(\EE)$ which acts properly on $\EE$
is virtually polycyclic.
 \item
A nonabelian free subgroup of $\Aff(\EE)$ cannot act properly
on $\EE$.
\item
The Euler characteristic 
$\chi(\Gamma\backslash\EE)$  (when defined)
of a complete affine manifold $\Gamma\backslash\EE$ 
must vanish (unless $\Gamma = \{1\}$ of course).
\item A complete affine manifold $\Gamma\backslash\EE$ has finitely
generated fundamental group $\Gamma$.
\end{itemize}
(If these conditions were met, one would have a satisfying
structure theory, similar to, but somewhat more involved, than
the Bieberbach structure theory for flat Riemannian manifolds.)

In \cite{Milnor}, Milnor provides abundant ``evidence'' 
for this ``conjecture''.
For example, the {\em infinitesimal version:\/}
Namely, let $G\subset\Aff(\EE)$ be a connected Lie group
which acts properly on $\EE$. Then $G$ must be an amenable Lie 
group, which simply means that it is a compact extension of a solvable
Lie group. (Equivalently, its Levi subgroup is compact.) 
Furthermore, he provides  a {\em converse:\/}
Milnor shows that every virtually polycyclic group
admits a proper affine action. (However, Milnor's actions
do {\em not\/} have compact quotient. Benoist~\cite{Benoist}
found finitely generated nilpotent groups which admit no
affine crystallographic action. 
Benoist's examples are $11$-dimensional.)

However convincing as his ``evidence'' is, Milnor still proposes
a possible way of constructing counterexamples:

\begin{quote}{\sl ``Start with a free discrete subgroup of $\Oo(2,1)$ and add translation components to obtain a group of affine transformations which  acts freely. 
However it seems difficult to decide whether the resulting group action is properly discontinuous.''}
\end{quote}
\noindent
This is clearly a geometric problem:
As Schottky showed in 1907,  
free groups act properly by isometries
on hyperbolic $3$-space, and 
hence by diffeomorphisms of $\EE^3$.
These actions are  {\em not\/} affine.

One might try to construct a proper affine action of a free group
by a construction like Schottky's.  
Recall that a {\em Schottky group of rank $g$\/} 
is defined by a system of $g$ open half-spaces
$H_1, \dots, H_g$ and isometries $A_1,\dots, A_g$ such that
the $2g$ half-spaces 
$$
H_1, \dots, H_g, A_1(H_1^c), \dots A_g(H_g^c)
$$ 
are all disjoint (where $H^c$ denotes the {\em complement\/} 
of the closure $\bar{H}$ of $H$). 
The {\em slab\/} 
$$
\mathsf{Slab}_i :=  H_i^c \cap A_i(H_i)
$$
is a fundamental domain for the action of the cyclic group
$\langle A_i\rangle$.
The {\em ping-pong lemma\/} then asserts that the intersection of 
all the slabs
$$
\Delta := \mathsf{Slab}_1 \cap \dots \cap \mathsf{Slab}_g
$$
is a fundamental domain for the group 
$\Gamma := \langle A_1, \dots, A_g\rangle$.
Furthermore $\Gamma$ is freely generated by
$A_1, \dots, A_g$.
The basic idea is the following.
Let $B_i^+ := A_i(H_i^c)$ 
(respectively $B_i^- := H_i$)
denote the {\em attracting basin\/} for $A_i$ (respectively $A_i^{-1}$).
That is, $A_i$ maps all of $H_i^c$ to $B_i^+$ and
$A_i^{-1}$ maps all of $A_i(H_i)$ to $B_i^-$.
Let $w(a_1, \dots, a_g)$ be a reduced word
in abstract generators $a_1,\dots, a_g$, with initial letter
$a_i^{\pm}$. 
Then 
$$
w(A_1,\dots,A_g) (\Delta) \subset B_i^{\pm}.
$$
Since all the basins $B_i^{\pm}$ are disjoint,
$w(A_1,\dots,A_g)$ maps $\Delta$ off itself.
Therefore $w(A_1,\dots,A_g) \neq 1$. 

Freely acting discrete cyclic groups of affine transformations 
have fundamental domains which are {\em parallel slabs,\/} that is,
regions bounded by two parallel affine hyperplanes. 
One might try to combine such slabs to form ``affine Schottky
groups'', but immediately one sees this idea is doomed, 
if one uses parallel slabs for Schottky's construction:
parallel slabs have disjoint complements only if they are
parallel to each other, in which case the group is necessarily
cyclic anyway. From this viewpoint, a discrete group of affine transformations seems very unlikely to act properly.

\section{Margulis spacetimes}

In the early 1980's Margulis, while trying to prove
that a nonabelian free group can't act properly by
affine transformations, 
discovered that  discrete free groups of affine transformations can indeed act properly!

Around the same time, David Fried and I were also working
on these questions, and reduced Milnor's question in dimension
three to what seemed at the time to be one annoying case which we could not handle. 
Namely, we showed the following: Let $\EE$ be a three-dimensional
affine space and $\Gamma\subset\Aff(\EE)$. 
Suppose that $\Gamma$ acts properly on $\EE$.
Then either $\Gamma$ is polycyclic or the restriction of
the linear holonomy homomorphism
$$
\Gamma \xrightarrow{\L} \GL(\EE)
$$
discretely embeds  $\Gamma$ onto a subgroup of $\GL(\EE)$
conjugate to the orthogonal group $\oto$.
 
In particular the complete affine manifold $M^3 = \Gamma\backslash\EE$
is a {\em complete flat Lorentz $3$-manifold\/} after one passes
to a finite-index torsionfree subgroup of $\Gamma$ to ensure
that $\Gamma$ acts freely. In particular the restriction $\L|_\Gamma$
defines a free properly discrete isometric action of $\Gamma$ on the
hyperbolic plane $\Ht$ and the quotient $\Sigma^2 := \Ht/\L(\Gamma)$
is a complete hyperbolic surface with a homotopy equivalence
$$
M^3 := \Gamma\backslash\EE \simeq \Ht/\L(\Gamma) =: \Sigma^2.
$$
Already this excludes the case when $M^3$ is compact,
since $\Gamma$ is the fundamental group of a closed
aspherical $3$-manifold (and has cohomological dimension $3$)
and the fundamental group of a hyperbolic surface (and has
cohomological dimension $\le 2$). 
This is a crucial step in the proof of 
Conjecture~\ref{conj:AuslanderConjecture}
in dimension $3$.

That the hyperbolic surface $\Sigma^2$ is 
{\em noncompact\/}
is a much deeper result due to Geoffrey Mess~\cite{Mess}.
Later proofs and a generalization have been found by
Goldman-Margulis~\cite{GoldmanMargulis} and Labourie~\cite{Labourie}. 
(Compare the discussion in  \S\ref{sec:Deformations}.)
Since the fundamental group of a noncompact surface
is free, $\Gamma$ is a free group. 
Furthermore $\L|_\Gamma$ embeds  $\Gamma$ 
as a free discrete group of isometries of hyperbolic space.
Thus Milnor's suggestion is the {\em only\/} way to construct
nonsolvable examples {\em in dimension three.\/}

 \subsection{Affine boosts and crooked planes}\label{sec:affineboosts}
Since $\L$ embeds $\Gamma_0$ as the fundamental group
of a hyperbolic surface, $\L(\gamma)$ is elliptic only if $\gamma = 1$.
Thus, if $\gamma\neq 1$,  then $\L(\gamma)$ is either hyperbolic or parabolic. Furthermore $\L(\gamma)$ is hyperbolic for most $\gamma\in\Gamma_0$.

When  $\L(\gamma)$ is hyperbolic, 
$\gamma$  is an {\em affine boost,\/} that is,
it has the form
\begin{equation}\label{eq:AffineBoost}
\gamma = \bmatrix e^{\ell(\gamma)} & 0 & 0 \\ 0 & 1 & 0 \\  
0 & 0 & e^{-\ell(\gamma)} \endbmatrix
\bmatrix 0 \\ \alpha(\gamma) \\  0 \endbmatrix
\end{equation}
in a suitable coordinate system.
(Here the $3\times 3$ matrix represents the linear part, 
and the column $3$-vector represents the translational part.)
$\gamma$ leaves invariant
a unique (spacelike) line $C_\gamma$ (the second coordinate line in 
\eqref{eq:AffineBoost}. 
Its image in $\Eto/\Gamma$ is a {\em closed geodesic\/} 
$C_\gamma/\langle\gamma\rangle$.
Just as for hyperbolic surfaces, most loops in $M^3$ are
freely homotopic to  such closed geodesics.
(For a more direct relationship between the dynamics of 
the geodesic flows on $\Sigma^2$ and $M^3$,
compare Goldman-Labourie~\cite{GoldmanLabourie}).

Margulis observed that $C_\gamma$ 
inherits a natural orientation and metric,
arising from an orientation on $\EE$, as follows.
Choose  repelling and attracting eigenvectors
$\L(\gamma)^\pm$ for $\L(\gamma)$ respectively; 
choose them so they lie in the same component of the nullcone. 
Then the orientation and metric on $C_\gamma$ is determined
by a choice of nonzero vector 
$\L(\gamma)^0$  spanning $\Fix\big(\L(\gamma)\big)$.
this vector is uniquely specified by requiring that:
\begin{itemize}
\item   $\L(\gamma)^0\cdot \L(\gamma)^0\ =\ 1$;
\item $\big( \L(\gamma)^0,  \L(\gamma)^-,
 \L(\gamma)^+\big)$ is a positively oriented basis.
\end{itemize}
The restriction of $\gamma$ to  $C_\gamma$ is
a translation by displacement $\alpha(\gamma)$
with respect to this natural orientation and metric.

Compare this to the more familiar {\em geodesic length function\/}
$\ell(\gamma)$ associated to a class $\gamma$ of closed curves
on the hyperbolic surface $\Sigma$.
The linear part $\L(\gamma)$ acts by {\em transvection\/}
along a geodesic $c_{\L(\gamma)}\subset\Ht$. 
The quantity $\ell(\gamma) > 0$
measures how far $\L(\gamma)$ moves
points of $c_{\L(\gamma)}.$ 

This pair of quantities 
$$
\big(\ell(\gamma),\alpha(\gamma)\big) \in
\R_+\times \R
$$
is a complete invariant of the isometry type of the
{\em flat Lorentz cylinder\/}
$\EE/\langle\gamma\rangle$.
The absolute value $\vert\alpha(\gamma)\vert$
is the length of the unique primitive closed geodesic in
$\EE/\langle\gamma\rangle$.

A fundamental domain is the parallel slab
$$
(\Pi_{C_\gamma})^{-1} 
\left( p_0\ +\  [0,\alpha(\gamma)]\ \gamma^0\right)
$$
where
$$
\EE\xrightarrow{\Pi_{C_\gamma}} C_\gamma
$$ 
denotes orthogonal projection onto 
$$
C_\gamma\ =\   p_0 + \R \gamma^0.
$$ As noted above, however, parallel slabs can't be combined to form
fundamental domains for Schottky groups, since their complementary
half-spaces are rarely disjoint.

In retrospect this is believable, since these fundamental
domains are fashioned from the dynamics of the 
translational part (using the projection $\Pi_{C_\gamma}$).
While the effect of the translational part is properness,
the dynamical behavior affecting most points is influenced
by the {\em linear part:\/} While points on $C_\gamma$
are displaced by $\gamma$ at a polynomial rate,
all other points move at an exponential rate. 

Furthermore, parallel slabs are less robust than slabs in $\Ht$:
while small perturbations of one boundary component extend
to fundamental domains, this is no longer true for parallel slabs.
Thus one might look for other types of fundamental
domains better adapted to the exponential growth
dynamics given by the linear holonomy $\L(\gamma)$.

Todd Drumm, in his 1990 Maryland thesis~\cite{Drumm_thesis}, 
defined more flexible polyhedral surfaces,  
which can be combined
to form fundamental domains for {\em Schottky groups\/}
of $3$-dimensional affine transformations.
A {\em crooked plane\/} is a PL surface in $\EE$, 
separating $\EE$ into two {\em crooked halfspaces.\/}
The complement of two disjoint crooked halfspace is a
{\em crooked slab,\/} which forms a fundamental domain for a cyclic
group generated by an affine boost. 
Drumm proved the remarkable theorem that if $S_1,\dots,S_g$
are crooked slabs whose complements have disjoint interiors,
then given any collection of affine boosts $\gamma_i$ with 
$S_i$ as fundamental domain, then the intersection
$S_1 \cap \dots\cap S_g$ is a fundamental domain for 
$\langle \gamma_1, \dots, \gamma_g\rangle$ acting on 
{\em all\/} of $\EE$. 

Modeling a crooked fundamental domain for 
$\Gamma$  acting on $\EE$ on a fundamental polygon
for $\Gamma_0$ acting on $\Ht$, Drumm proved the following
sharp result:

\begin{thm*}[Drumm~\cite{Drumm_thesis,Drumm_jdg}]
Every {\em noncocompact\/} torsionfree Fuchsian group $\Gamma_0$ 
admits a proper affine deformation $\Gamma$ whose quotient
is a solid handlebody.
\end{thm*}
\noindent
(Compare also \cite{DrummGoldman1,CharetteGoldman}.)

\begin{figure}[h]
\centering
\mbox{\subfigure
{\epsfxsize=2.4in \epsfbox{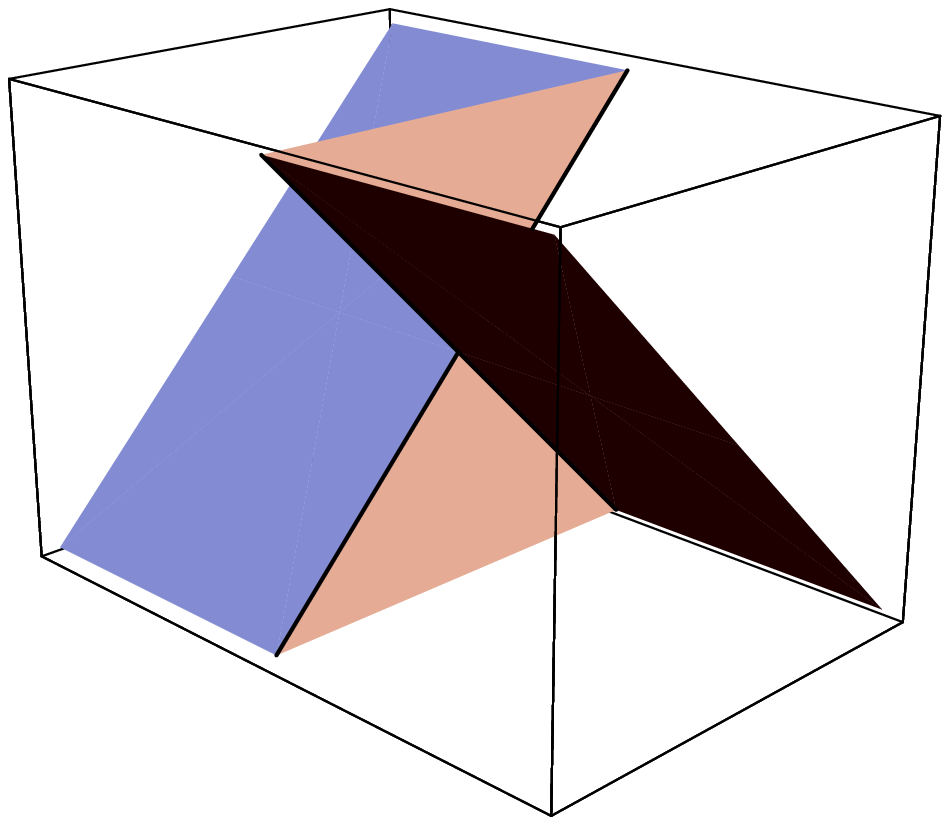}}
}
\mbox{\subfigure
{\epsfxsize=2.4in \epsfbox{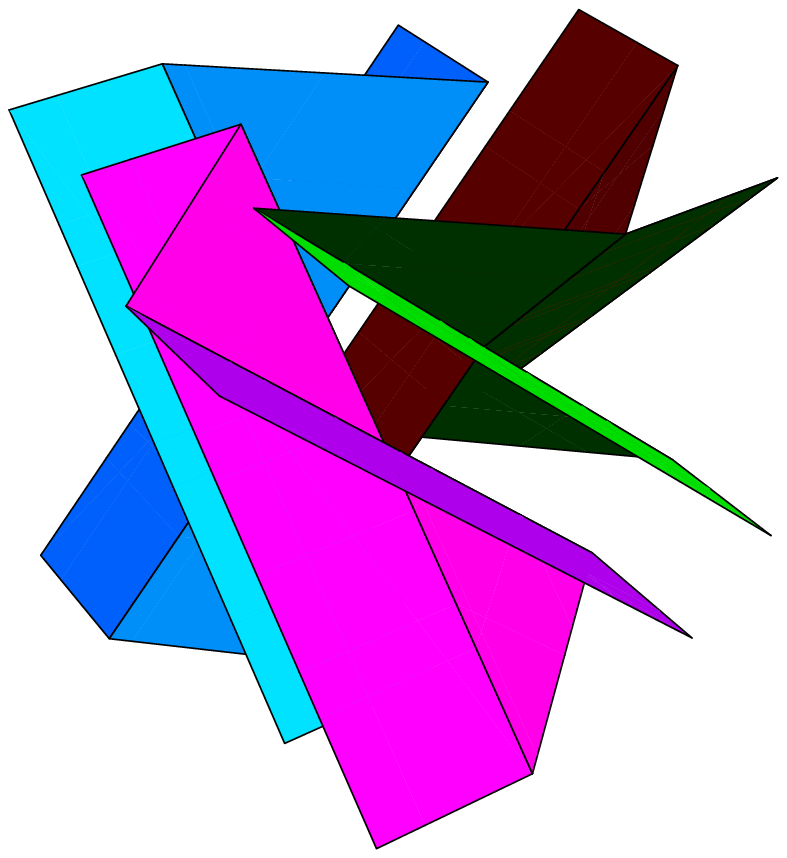}} 
}
\caption{A crooked plane, and a family of three pairwise disjoint crooked planes}
\end{figure}

\newpage
\subsection{Marked length spectra}\label{sec:MarkedLengthSpectra}

We now combine the geodesic length function $\ell(\gamma)$
describing the geometry of the hyperbolic surface $\Sigma$
with the Margulis invariant $\alpha(\gamma)$  describing
the Lorentzian geometry of the flat affine $3$-manifold $M$.

As noted by Margulis, $\alpha(\gamma) = \alpha(\gamma^{-1})$,
and more generally
$$
\alpha(\gamma^n)\ =\  \vert n\vert \alpha(\gamma).
$$
The invariant $\ell$ satisfies the same homogeneity condition,
and therefore
$$
\frac{\alpha(\gamma^n)}{\ell(\gamma^n)}\ =\ 
\frac{\alpha(\gamma)}{\ell(\gamma)}
$$
is constant along hyperbolic cyclic subgroups.
Hyperbolic cyclic subgroups correspond to 
periodic orbits of the geodesic flow $\phi$ on the
unit tangent bundle $U\Sigma$.
Periodic orbits, in turn, 
define $\phi$-invariant probability measures on $U\Sigma$.  
Goldman-Labourie-Margulis~\cite{GLM}
prove that, for any affine deformation,
this function extends to a continuous function $\Upsilon_\Gamma$
on the space $\CS$ of
$\phi$-invariant probability measures on $U\Sigma$. 
Furthermore when $\Gamma_0$ is convex cocompact 
(that is, contains no parabolic elements), 
then the affine deformation $\Gamma$ acts properly
if and only if $\Upsilon_\Gamma$ never vanishes.
Since $\CS$ is connected,
nonvanishing implies either all $\Upsilon_\Gamma(\mu) > 0$
or all $\Upsilon_\Gamma(\mu) < 0$.
From this follows Margulis's {\em Opposite Sign Lemma,\/}
first proved in \cite{Margulis1,Margulis2} and extended
to groups with parabolics by Charette and Drumm~\cite{CharetteDrumm}:

\begin{thm*}[Margulis]
If $\Gamma$ acts properly, then all of the numbers $\alpha(\gamma)$
have the same sign.
\end{thm*}

\noindent
For an excellent treatment of the original proof of this fact, 
see the survey article of Abels~\cite{Abels}.

\subsection{Deformations of hyperbolic surfaces}\label{sec:Deformations}

The Margulis invariant may be interpreted in terms of deformations
of hyperbolic structures as follows~\cite{GoldmanMargulis, Goldman_MargInv}).

Suppose $\Gamma_0$ is a Fuchsian group with quotient
hyperbolic surface $\Sigma_0 = \Gamma_0\backslash\Ht$.
Let $\gad$ be the $\Gamma_0$-module defined by the 
adjoint representation applied to the embedding
$\Gamma_0\hookrightarrow\oto$.
The coefficient module $\gad$ corresponds to the Lie algebra
of {\em right-invariant\/} vector fields on $\oto$ with the action
of $\oto$ by left-multiplication. Geometrically these vector fields
correspond to the infinitesimal isometries of $\Ht$.

A family of hyperbolic surfaces
$\Sigma_t$  smoothly varying with respect to a parameter $t$
determines an  {\em infinitesimal deformation,\/} which is 
a cohomology class  $[u]\in H^1(\Gamma_0,\gad)$,
The cohomology group $H^1(\Gamma_0,\gad)$ corresponds to 
{\em infinitesimal deformations\/} of the hyperbolic surface $\Sigma_0$.
In particular the tangent vector to the path $\Sigma_t$
of marked hyperbolic structures smoothly varying with respect to a parameter $t$
defines a cohomology class 
$$
[u]\in H^1(\Gamma_0,\gad).
$$

The same cohomology group parametrizes affine deformations.
The translational part $u$ of a linear representations of $\Gamma_0$ 
is a cocycle of the group $\Gamma_0$ taking values 
in the corresponding $\Gamma_0$-module $\V$. 
Moreover two cocycles define affine deformations which are
conjugate by a translation if and only if their translational parts
are cohomologous cocycles. Therefore translational
conjugacy classes of affine deformations form  the
cohomology group $\ho$. Inside $\ho$ is the subset $\Proper$
corresponding to {\em proper affine deformations.}

The adjoint representation $\Ad$ of $\oto$ identifies with the
orthogonal representation of $\oto$ on $\V$. 
Therefore  the cohomology group $\ho$ consisting of translational conjugacy
classes of affine deformations of $\Gamma_0$ 
can be identified %
with the 
cohomology group $H^1(\Gamma_0,\gad)$ corresponding to 
infinitesimal deformations of $\Sigma_0$. 

\begin{thm*}
Suppose $u\in Z^1(\Gamma_0,\gad)$ defines an
{\em infinitesimal deformation\/} tangent to
a smooth deformation
$\Sigma_t$ of $\Sigma$.
\begin{itemize}
\item The marked length spectrum $\ell_t$ of $\Sigma_t$
varies smoothly with $t$.
\item
Margulis's invariant $\alpha_u(\gamma)$ represents the derivative
\begin{equation*}
\frac{d}{dt}\bigg|_{t=0} \ell_t(\gamma)
\end{equation*}
\item (Opposite Sign Lemma)
If $[u]\in \Proper$, then all closed geodesics lengthen
(or shorten) under the deformation $\Sigma_t$. 
\end{itemize}
\end{thm*}
\noindent
Since closed hyperbolic surfaces do not support deformations
in which {\em every\/} closed  geodesic shortens,
such deformations only exist when $\Sigma_0$ is noncompact.
This leads to a new proof~\cite{GoldmanMargulis} of 
Mess's theorem that $\Sigma_0$ is not compact.
(For another, somewhat similar proof, which generalizes to higher dimensions,
see Labourie~\cite{Labourie}.)

The tangent bundle $TG$ of any Lie group $G$ has a natural structure
as a Lie group, where the fibration $TG\xrightarrow{\Pi} G$ is a
homomorphism of Lie groups, and the tangent spaces 
$$
T_xG = \Pi^{-1}(x) \subset TG
$$ 
are vector groups.  The deformations of a representation
$\Gamma_0\xrightarrow{\rho_0}G$ correspond to representations
$\Gamma_0\xrightarrow{\rho}TG$ such that $\Pi\circ\rho = \rho_0$.
In our case, affine deformations of $\Gamma_0\hookrightarrow\oto$ 
correspond to representations in the tangent bundle $T\oto$.
When $G$ is the group $\mathsf{G}(\R)$ of $\R$-points
of an algebraic group $\mathsf{G}$ defined over $\R$, then 
$$
TG\ \cong\  \mathsf{G}(\R[\epsilon])
$$ 
where $\epsilon$ is an indeterminate with $\epsilon^2 = 0$.
(Compare~\cite{Goldman_MargInv}.)
This is reminiscent of the classical theory of quasi-Fuchsian deformations
of Fuchsian groups, where one deforms a Fuchsian subgroup of 
$\SLtR$
in $$
\mathsf{SL}(2,\C) = \mathsf{SL}(2,\R[i])
$$ 
where $i^2 = -1$.

\subsection{Classification}

In light of Drumm's theorem, classifying  Margulis spacetimes $M^3$
begins with 
the classification of hyperbolic structures $\Sigma^2$.
Thus the deformation space of Margulis spacetimes maps to
the Fricke space $\Fricke$ of marked hyperbolic structures on
the underlying topology of $\Sigma$.

The main result of \cite{GLM} 
is that the positivity (or negativity) of  $\Upsilon_\Gamma$ on 
on $\CS$ is necessary and sufficient for
properness of $\Gamma$. (For simplicity we restrict ourselves
to the case when $\L(\Gamma)$ contains no parabolics ---
that is, when $\Gamma_0$ is convex cocompact.) 
Thus the proper affine deformation space $\Proper$ identifies
with %
the open convex cone in $\ho$ defined by the linear functionals
$\Upsilon_\mu$, for $\mu$ in the compact space $\CS$.
These 
give %
uncountably many linear conditions on $\ho$,
one for each $\mu\in\CS$.
Since the invariant probability measures arising from periodic
orbits are dense in $\CS$, the cone $\Proper$ is the interior of half-spaces
defined by the countable set of functional $\Upsilon_\gamma$,
where $\gamma\in\Gamma_0$. 

The zero level sets $\Upsilon_\gamma^{-1}(0)$ %
correspond to affine deformations
where $\gamma$ does not act freely. %
Therefore $\Proper$ defines a component of the subset
of $\ho$ corresponding to affine deformations 
which are {\em free\/} actions.

\begin{figure}[H]
\centering
\mbox{\subfigure
[Three-holed sphere]
{\epsfxsize=2.4in \epsfbox{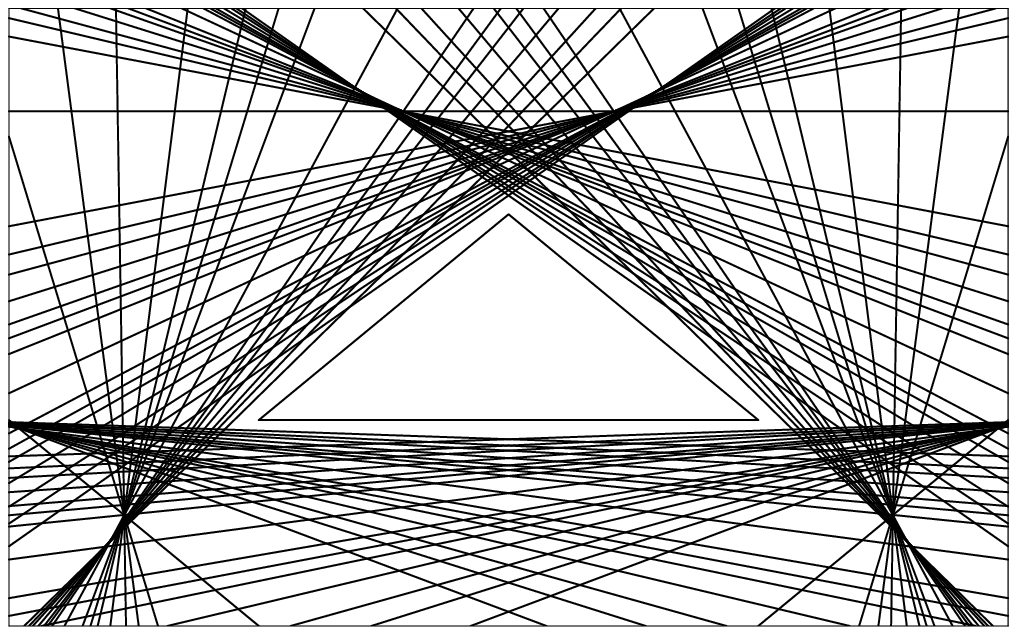}}
}
\mbox{\subfigure[Two-holed cross-surface $\R\mathsf{P}^2$]
{\epsfxsize=2.4in \epsfbox{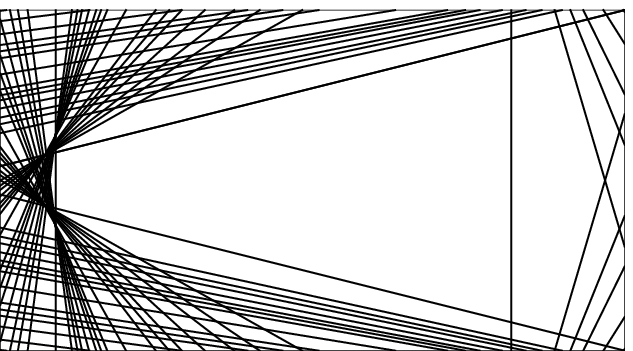}}
}
\caption{Finite-sided deformation spaces for surfaces with $\chi(\Sigma)$}
\label{fig:FinDefs}
\end{figure}

\begin{figure}[H]
\centering
\mbox{\subfigure[One-holed torus]
{\epsfxsize=2.4in \epsfbox{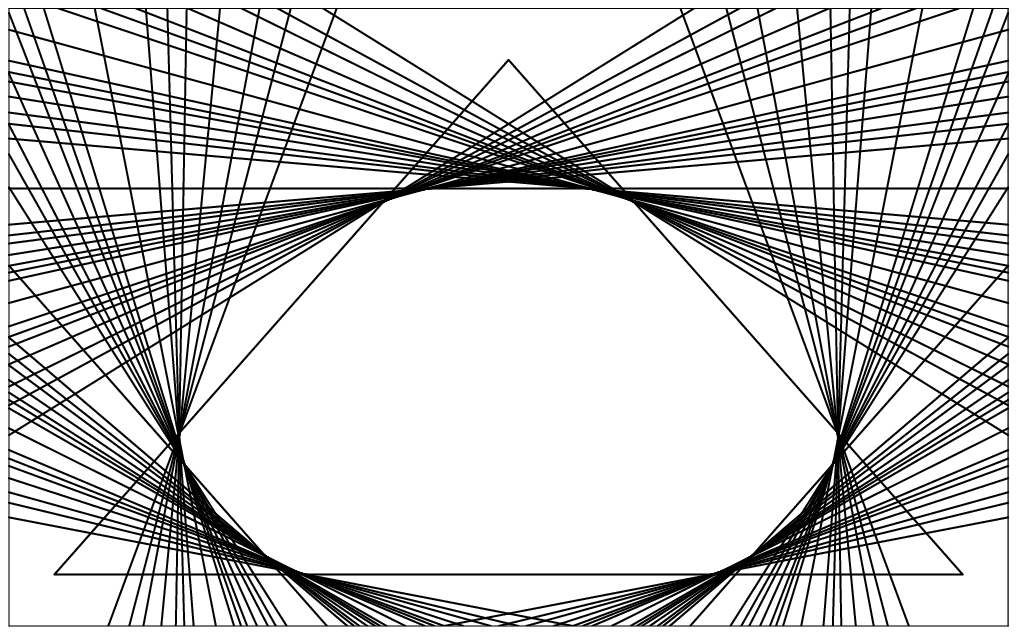}}
}
\mbox{\subfigure[One-holed Klein bottle]
{\epsfxsize=2.4in \epsfbox{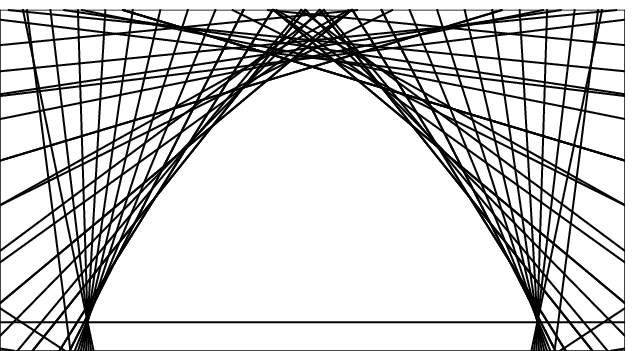}}
}
\caption{Infinite-sided deformation spaces for surfaces with $\chi(\Sigma)$}
\label{fig:InfDefs}
\end{figure}

Actually, one may go further.  An argument inspired by
Thurston~\cite{Thurston_stretch}, reduces the consideration to only
those measures arising from {\em multicurves,\/} that is, unions of
disjoint {\em simple\/} closed curves. These measures (after scaling)
are dense in the {\em Thurston cone\/} $\ml$ of {\em measured geodesic laminations\/} 
on $\Sigma$.  One sees the combinatorial structure of
the Thurston cone replicated on the boundary of $\Proper\subset\ho$.
(Compare Figures~\ref{fig:FinDefs} and \ref{fig:InfDefs}.)


Two particular cases are  notable. When $\Sigma$ is a
$3$-holed sphere or a $2$-holed cross-surface (real projective
plane), then the Thurston cone degenerates to a finite-sided
polyhedral cone. In particular properness is characterized by  $3$ 
Margulis functionals for the $3$-holed sphere,
and $4$ for the $2$-holed cross-surface. Thus the deformation space
of equivalence classes of proper affine deformations is either
a cone on a triangle or a convex quadrilateral, respectively.

When $\Sigma$ is a $3$-holed sphere, 
these functionals correspond to the three components of
$\partial\Sigma$. 
The halfspaces defined by the correponding three Margulis functionals
cut off the deformation space (which is a polyhedral cone with $3$ faces).
The Margulis functionals for the other curves define halfspaces which
strictly contain this cone. 

When $\Sigma$ is a $2$-holed cross-surface these functionals correspond
to the two components of $\partial\Sigma$ and the two
orientation-reversing simple closed curves in the interior of
$\Sigma$. The four Margulis functionals describe a polyhedral cone
with $4$ faces.  All other closed curves on $\Sigma$ define halfspace
strictly containing this cone.  

In both cases, an ideal triangulation for $\Sigma$ models a crooked
fundamental domain for $M$, and $\Gamma$ is an affine Schottky group,
and $M$ is an open solid handlebody of genus $2$
(Charette-Drumm-Goldman~\cite{CDG,CDG1,CDG2}).  Fig.~\ref{fig:FinDefs}
depicts these finite-sided deformation spaces.

For the other surfaces where $\pi_1(\Sigma)$ is free of rank two
(equivalently $\chi(\Sigma) = -1$), infinitely many functionals
$\Upsilon_\mu$ are needed to define the deformation space, which
necessarily has infinitely many sides. In these cases $M^3$
admits crooked fundamental domains corresponding to ideal triangulations
of $\Sigma$, although unlike the preceding cases there is no single
ideal triangulation which works for all proper affine deformations.
Once again $M^3$ is a genus two handlebody.
Fig.~\ref{fig:InfDefs} depicts these infinite-sided deformation spaces.

\subsection{An arithmetic example}

These examples are everywhere. As often happens in mathematics,
finding the first example of generic behavior can be quite difficult.
However, once the basic phenomena are recognized, examples of this
generic behavior abound.  The following example, taken from
\cite{CDG1}, shows how a proper affine deformation sits inside the
symplectic group $\mathsf{Sp}(4,\Z)$ .

Begin with a $2$-dimensional vector space $L_0$ over $\R$ with
the group of linear automorphisms $\GL(L_0)$. 
Let $\V$ denote the vector space of {\em symmetric bilinear maps\/}
$L_0\times L_0 \xrightarrow{b}\R$ with the induced action of $\GL(L)$.
Identifying $\V$ with symmetric $2\times 2$ real matrices,
the negative of the determinant defines an invariant Lorentzian
inner product on $\V$. In particular this defines a local  embedding
$\GL(L_0) \longrightarrow \oto$.

Let $\Li := L^*$ denote the vector space dual to $L_0$ and $\W := L_0 \oplus \Li$
the direct sum. Then $\W$ admits a unique symplectic structure $\omega$
such that $L_0$ and $\Li$ are Lagrangian subspaces and the restriction of $\omega$ 
to $L_0\times \Li$ is the duality pairing. 
Let $\mathsf{Sp}(4,\R)$ denote the group of linear symplectomorphisms of
$(\W,\omega)$. It acts naturally on the homogeneous space $\mathscr{L}(\W,\omega)$
of Lagrangian $2$-planes $L$ in $(\W,\omega)$. 

The Minkowski space $\EE$ associated to $\V$ consists of Lagrangians
$L\in\mathscr{L}(\W,\omega)$ which are transverse to $\Li$. This is a
torsor for the Lorentzian vector space $\V$ as follows.  $\V$ consists
of symmetric bilinear forms on $L_0$, and these 
can be identified with {\em self-adjoint\/} linear maps 
$$
L_0 \xrightarrow{f}\Li\cong (\L_0)^*.
$$
A $2$-dimensional linear subspace of $\V$ which is transverse to $\Li$ is the graph
$L := \mathsf{graph}(f)$ of a linear map $L_0\xrightarrow{f}\Li$. 
Moreover $L$ is Lagrangian
if and only if $f$ is self-adjoint. 
Furthermore, since $\V$ is a vector space, it acts simply transitively on the space $\EE$
of such graphs by addition. In terms of $4\times 4$ symplectic matrices ($2\times 2$ block
matrices using the decomposition $\W = L_0 \oplus \Li$), 
these translations correspond to {\em shears:\/} 
\begin{equation}\label{eq:shears}
\bmatrix I_2 & f_2 \\ 0 & I_2 \endbmatrix
\end{equation}
where the corresponding symmetric $2\times 2$ matrix corresponding
to $f$ is denoted $f_2$.  The corresponding subgroup of
$\mathsf{Sp}(4,\R)$ consists of linear symplectomorphisms of
$(\W,\omega)$ which preserve $\Li$, and act identically both on $\Li$
and on its quotient $\W/\Li$.

As the translations of $\EE$ are represented by shears 
in block upper-triangular form~\eqref{eq:shears},
the linear isometries are represented by the block diagonal matrices
arising from $\mathsf{SL}(L_0)$. 
More generally, the {\em Lorentz similarities\/} of $\EE$ correspond
to $\mathsf{GL}(L_0)$ as follows. A linear automorphism $L_0\xrightarrow{g}L_0$
induces a linear symplectomorphism 
$g\oplus (g^\dagger)^{-1}$ of $\W = L_0 \oplus \Li$:
$$
g\oplus (g^\dagger)^{-1} \;=\;
\bmatrix g & 0 \\ 0 & (g^\dagger)^{-1} \endbmatrix.
$$
These linear symplectomorphisms can be characterized as those which
preserve the Lagrangian $2$-planes $L_0$ and $\Li$.
Furthermore $g$ induces an {\em isometry\/} of $\EE$ with the flat Lorentzian structure
if and only if $\mathsf{Det}(g) = \pm 1$. 

Here is our example.  The level two congruence subgroup $\Gamma_0$ is
the subgroup of $\mathsf{GL}(2,\Z)$ generated by
$$
\bmatrix  -1 & -2 \\
0 & -1 \endbmatrix ,\;
\bmatrix
-1 &  0 \\
 2 & -1 \endbmatrix
$$
and the corresponding hyperbolic surface is a triply punctured sphere.
For $i=1,2,3$ choose three positive integers $\mu_1,\mu_2,\mu_3$
(the coordinates of the translational parts).
Then the subgroup $\Gamma$ of $\mathsf{Sp}(4,\Z)$ generated by
\begin{equation*}
\bmatrix-1
 & -2 & \mu_1 + \mu_2 -\mu_3  & 0 \\
 0 & -1 & 2\mu_1 & -\mu_1 \\
 0 & 0 & -1   &   0 \\
 0 & 0 &  2   &   -1
\endbmatrix ,\;
\bmatrix
 -1 &  0 & -\mu_2 & -2\mu_2    \\
  2 & -1 &     0 & 0 \\
  0 & 0  &      -1 & -2   \\
  0 & 0  &      0 & -1
 \endbmatrix
 \end{equation*}
defines a affine deformation of a $\Gamma_0$.

\begin{figure}[ht]
\centerline{\epsfxsize=5.0in \epsfbox{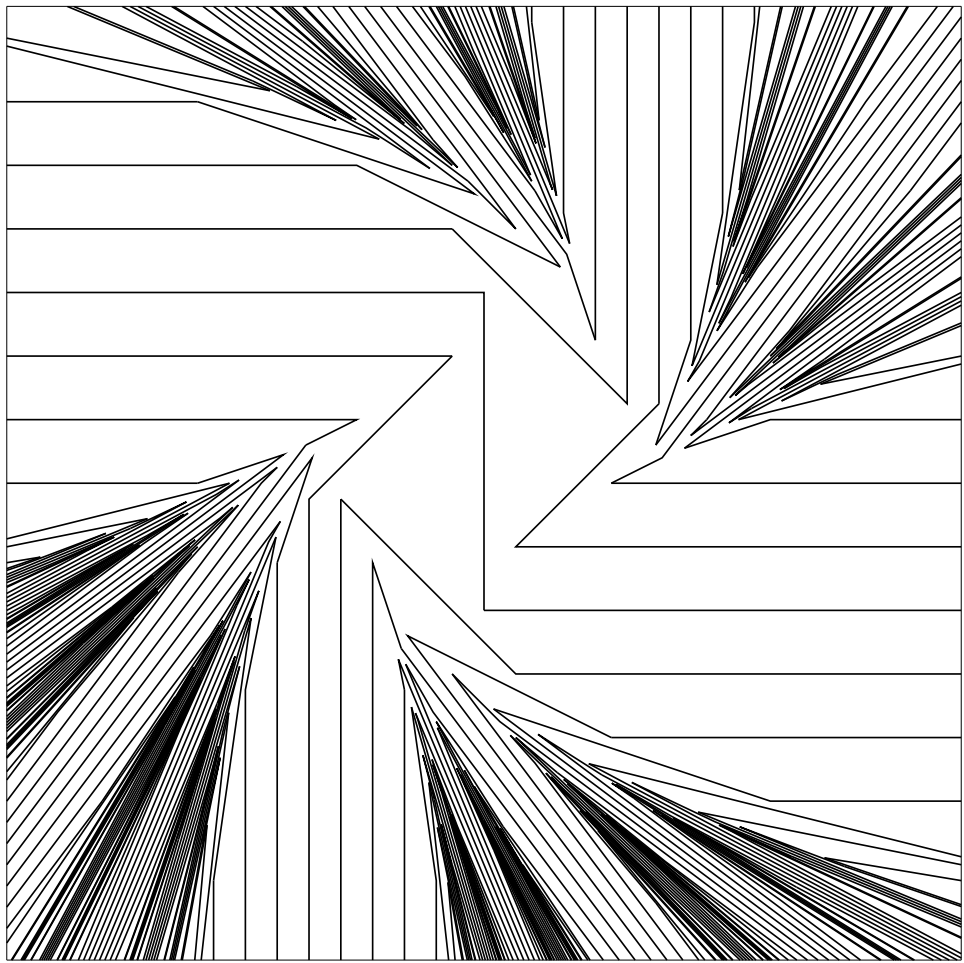}}
\caption{
A proper affine deformation of level two congruence subgroup
of $\mathsf{SL}(2,\Z)$.}
\label{fig:modularzigzags}
\end{figure}

By the main result of Charette-Drumm-Goldman~\cite{CDG}, this affine
deformation is proper with a fundamental polyhedron bounded by crooked
planes.  The quotient $3$-manifold $M^3 = \Gamma\backslash\EE$ is
homeomorphic to a genus two handlebody.  Fig.~\ref{fig:modularzigzags}
depicts the intersections of crooked fundamental domains for this
group (when $\mu_1=\mu_2=\mu_3 = 1$) with a spacelike plane.  Note the
parallel line segments cutting off fundamental domains for the cusps
of $\Sigma$; the parallelism results from the parabolicity of the
holonomy around the cusps.

\end{document}